\documentclass[10pt,twocolumn,letterpaper]{article}

\usepackage{cvpr}
\usepackage{times}
\usepackage{epsfig}
\usepackage{graphicx}
\usepackage{amsmath}
\usepackage{amssymb}

\usepackage{caption}
\usepackage{subcaption}

\usepackage[breaklinks=true,bookmarks=false]{hyperref}

\cvprfinalcopy 

\def\cvprPaperID{****} 

\begin{document}

\title{Fast uncertainty quantification of reservoir simulation with variational U-Net}

\author{Larry Jin, \quad Hannah Lu,\quad Gege Wen\\
Energy Resources Engineering\\
Stanford University\\
{\tt\small \{zjin, hannahlu, gegewen\}@stanford.edu}}

\maketitle

\begin{abstract}
Quantification of uncertainty in production/injection forecasting is an important aspect of reservoir simulation studies. Conventional approaches include intrusive Galerkin-based methods (e.g., generalized polynomial chaos (gPC)  and stochastic collocation (SC) methods) and non-intrusive Monte Carlo (MC) based methods. Nevertheless, the quantification is conducted in reformulations of the underlying stochastic PDEs with fixed well controls. If one wants to take various well control plans into account, expensive computations need to be repeated for each well design independently. In this project, we take advantages of the \textbf{equation-free} spirit of convolutional neural network (CNN) to overcome this challenge and thus achieve the flexibility of efficient uncertainty quantification with various well controls. We are interested in the development of surrogate models for uncertainty quantification and propagation in reservoir simulations using a deep convolutional encoder-decoder network as an analogue to the  image-to-image regression tasks in computer science. First, a U-Net architecture is applied to replace conventional expensive deterministic PDE solver. Then we adopt the idea from shape-guided image generation using variational U-Net and design a new variational U-Net architecture for \textbf{``control-guided"} reservoir simulation. Backward propagation is learned in the network to extract the hidden physical quantities and then predict the future production by the learned forward propagation using the hidden variable with various well controls. Comparisons in computational efficiency are made between our proposed CNN approach and conventional MC approach. Significant improvements in computational speed with reasonable accuracy loss are observed in the numerical tests. 
\end{abstract}

\section{Introduction}
Reservoir simulation is an important tool for a wide range of engineering problems including nuclear waste management, oil and gas production, as well as carbon capture and storage (CCS). It helps us understanding these subsurface multiphase flow processes by solving spatially and temporally discretized mass and energy balance equations. However, the veracity of reservoir simulation results is uncertain due to the fact that the model parameters are subject to the uncertainties of various sources. Sources of parameter uncertainties include errors in data used for model parameterization, interpretive errors, and most importantly, heterogeneity of the subsurface under consideration. 

The goal of uncertainty quantification (UQ) in reservoir simulation is to evaluate the impact of parameter uncertainties on the model outputs. Conventional methods for UQ are mostly Monte Carlo (MC) based sampling \cite{robert2013monte}, which could induce huge computational costs in repetitive multi-scale/multi-physics deterministic simulations \cite{giles2015multilevel}. Moreover, well location, as another important control in reservoir simulation, needs to be considered in the evaluation of production/injection proposals. MC based method requires enough sampling in each fixed well location, which involves time-consuming computations. To mitigate the challenges in computation cost, and allow efficient quantification of uncertainty in flow dynamics with high dimensional random media, we introduce a CNN-based surrogate modeling tool (i.e., a variational U-net) in this project.

CNN-based generative adversarial network (GAN), autoregressive model, and auto-encoder (AE) are efficient tools to conduct image-to-image regressions. The recent applications of CNN-based models in reservoir simulation have achieved encouraging results \cite{zhu2018bayesian, mo2018deep}. In those models, the input reservoir properties (e.g., permeability) and output flow states (e.g., pressure and saturation) are treated as images.The CNN-based image-to-image regression models are therefore used as surrogate models for UQ tasks. However, the aforementioned frameworks are not able to incorporate information from controlling parameters such as bottom-hole pressure (BHP) or flow rates, which limit their potential in the UQ of reservoir models where control plays a distinct role. 

A variational U-net has been proposed by \cite{esser2018variational} in the context of image and video generation. A conditional U-net \cite{ronneberger2015u} for shape-guided image generation, conditioned on the output of a variational auto-encoder for appearance is presented, where the model learns to infer appearance from the queries and synthesizes images with specific appearances in different poses (i.e., shapes).

In this project, we propose to adopt both the U-Net and the variational U-Net as surrogate models for efficient UQ tasks in reservoir simulations. The shape-guided image generation process with variational U-Net is analogous to \textbf{``control-guided" reservoir simulation}, and can be intuitively adopted for our tasks. Some previous studies related to this project are summarized in section 2. In section 3, we set up the problem in a rigorous mathematical frame and explain the dataset in details. In section 4, we propose a new variational U-Net architecture for control-guided reservoir simulation. Motivations and detailed explorations are also illustrated in subsections. Numerical results and comparison with conventional MC methods are presented in section 5. Finally, we draw a conclusion of this project and discuss some future research directions in section 6.

\section{Related Work}
The subsurface is such a complex system involving interactive physics in multiple scales that there is no accurate deterministic model. Thus probabilistic models have been explored to account for the uncertainties from model errors, model parametrization, heterogeneity of the environment and various geometries of the boundary, giving rise to extensive research interests in uncertainty quantification in reservoir engineering. A dominant strategy for such problems is to solve the deterministic problem at a finite large number of realizations of the random inputs using Monte Carlo sampling~\cite{robert2013monte}. Variations of MC methods, including quasi-MC~\cite{kalla2005use}, Multi-Level-MC(MLMC)~\cite{giles2015multilevel,lu2016improved}, stratified MC~\cite{burhenne2011sampling}, and etc, are designed to make more efficient sampling. Intrusive methods, (e.g., moment equations~\cite{winter2003moment}, polynomial chaos~\cite{xiu2003modeling}, stochastic collocation~\cite{babuvska2007stochastic}, method of distributions~\cite{tartakovsky2016method}), as alternatives to MC simulations, have been widely studied in the past decades. Although efficient in some problems, intrusive methods, in general, require additional efforts in reformulation of the model and reconstruction of the deterministic solvers. In the meantime, it is well known that all intrusive methods suffer from curse of dimensionality (in random space). We refer to reviews on these topics in~\cite{xiu2010numerical}.

Recently, deep learning has been explored as a competitive methodology across fields such as fluid mechanics~\cite{kutz2017deep}, hydrology~\cite{marccais2017prospective,chan2018machine}, bioinformatics~\cite{min2017deep}, high energy physics~\cite{baldi2014searching} and others. In particular, \cite{zhu2018bayesian} adopted an end-to-end image-to image regression approach for surrogate modeling governed by stochastic PDEs with high-dimensional stochastic input in random porous media. In addition, the deep neural networks are set under a formal Bayesian formula  to enable the network to express its uncertainty on its predictions when using limited training data. \cite{zhu2018bayesian} and our project both study two-dimensional, single phase, steady-state flow through a random permeability field. However, the emphasis of \cite{zhu2018bayesian} lies in Bayesian deep learning dealing with high dimensional random inputs but with fixed single well control. In this work, we focus on uncertainty propagation with \textbf{various well controls}, which would allow existing well information to infer the incomplete/uncertain geological properties and evaluate potential well locations during simulation. In large-scale reservoir project like CCS, such quantification from simulations will be of significant importance and economic value in industry.

In CS231N class, we have learned two different approaches to image generation in the context of deep learning: Variational Auto-Encoder (VAE)~\cite{kingma2013auto} and Generative Adversarial Networks (GAN)~\cite{goodfellow2014generative}. In \cite{esser2018variational}, a conditional U-Net for shape-guided image generation, conditioned on the output of a VAE for appearance is presented. The separation between shape and appearance is carefully modelled and thus an explicit representation of the appearance, which can be combined with new poses, is obtained. Motivated by this work, we established the analogue between shape-guided image and control-guided saturation map. Similarly, the hidden appearance in image is analogous to the underlying permeability map. Under this framework, the variational U-Net is trained to learn the underlying Darcy's law. We take advantages of the \textbf{equation-free} nature of CNN, allowing the freedom of saturation map generation under various well controls. Efficient uncertainty quantification tasks can be conducted then to provide valuable evaluations of potential injection plans.

\section{Dataset}
In this project, we consider a carbon capture and storage (CCS) problem as the reservoir simulation model of interest, where the multiphase phase flow consists carbon dioxide (CO$_{2}$) as gas and water as fluid. CCS is an essential climate change mitigation technology to reduce the concentration of CO$_{2}$ in the atmosphere. CO$_{2}$ captured from concentrated sources, the atmosphere, or through bio-energy production, is compressed into a liquid and injected into deep geological formations for long term sequestration. \cite{IPCC}

\subsection{Problem statement}
Let $S$ be a image of CO$_2$ gas saturation map from a dataset $\mathcal S$, which consists of different saturation samples from different injection well locations in various geological formation. The geological properties of interest are characterized by the permeability field $\mathbf k$ and the injection well location is denoted by $\mathbf y$. Physically, the relationship between the three quantities $S,\mathbf y$ and $\mathbf k$ are characterized by the following Eq.~\ref{eq:1}:
\begin{equation}\label{eq:1}
\begin{aligned}
&\phi\partial_t S+\nabla\cdot \mathbf u =0\\
&\mathbf u(\mathbf x) = -\mathbf k(\mathbf x,\omega)\nabla \Phi(\mathbf x), \mathbf x\in\mathcal D,\\
&\nabla\cdot \mathbf u(\mathbf x) = f(\mathbf x),\mathbf x\in\mathcal D,\\
&\mathbf u(\mathbf x) \cdot \hat{\mathbf n}(\mathbf x)=0 , \mathbf x\in\partial D,\\
\end{aligned}
\end{equation} 
where $\phi$ is the porosity, $\mathbf u$ is the velocity field and $\Phi$ is the potential including the pressure and gravity effect. The domain of interest (after projection) is $\mathcal D = [0,l]\times[0,w]$, where $l$ is the length and $w$ is the width of the system. The uncertainties of the system are indicated as the random variable $\omega$ in $\mathbf k(\mathbf x, \omega)$. Without loss of generality, we consider a Gaussian permeability field $\mathbf k$ with correlation length $25$. In our simulation, the reservoir is rescaled in a $128\times 128$ mesh with $10$-meter layer depth.  The spacial variable $\mathbf x = (x_1,x_2)$ specifies the coordinates. $f(\mathbf x)$ describes the injection by
\begin{equation}\label{eq:2}
    f(\mathbf x) = \left\{
    \begin{aligned}
    &Q, &\mbox{if} \ \mathbf x= \mathbf y,\\
    &0, &\mbox{otherwise},
    \end{aligned}
    \right.
\end{equation}
where $Q$ is the injection rate.

We want to understand how the saturation maps $S$ are influenced by the permeability map $\mathbf k$ and well location map $\mathbf y$ without repeatedly solving the equation system~(\ref{eq:1}). For an arbitrary given control of $\mathbf y$, how can we predict the uncertainty propagation in the saturation map $S$ will look like by using convolutional neural network? This can be interpreted as very useful information in the evaluation of potential CCS site.

\subsection{Data generation} 
The data set is generated by our team and the data preparation workflow is described as follows: 1) generate Gaussian permeability maps $\mathbf k$; 2) prepare the random well location map $\mathbf y$; 3) feed $\mathbf k$ and $\mathbf y$ to a numerical simulator and obtain the simulation results of CO$_2$ gas saturation map $S$ and pressure map $P$; and 4) concatenate two sets of $\mathbf y$, $S$, and $P$ to make one data tuple.

\begin{figure}[!htp]
    \centering
    \includegraphics[scale=0.5]{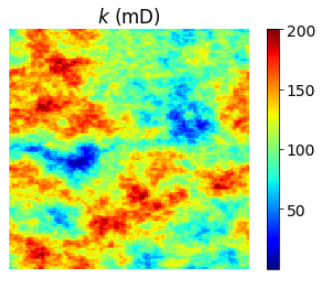}
    \caption{Example of the Gaussian permeability field $\mathbf k$}
    \label{fig:perm}
\end{figure}

In Figure~\ref{fig:perm}, we demonstrate an example of the Gaussian permeability fields $\mathbf k$ generated using the Stanford Geostatistical Modeling Software (SGeMS) \cite{sgems}. The permeability varies from 0.001 mD to 200 mD to mimic the sandstone reservoir that can potentially be used for CCS projects.

\begin{figure}[!htp]
  \centering
  \includegraphics[scale = 0.35]{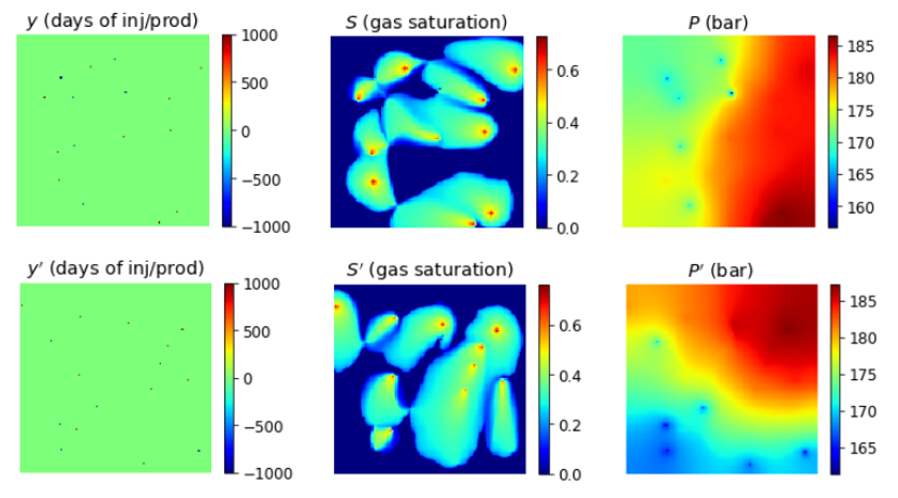}
  \caption{Example of a data tuple: $\mathbf y$, $S$, $P$, $\mathbf y'$, $S'$, and $P'$}
  \label{fig:data}
\end{figure}

The well location map $\mathbf{y}$ incorporates the injection and the production information. There are a total of 10 injection wells on each well location map and each injection well has a constant supercritical CO$_{2}$ flow rate of 75 tons per day. Each map contains 6 production wells with a constant bottom hole pressure of 150 bar. An example of $\mathbf{y}$ is shown in the first column of Figure~\ref{fig:data}. Note that the positive numbers of days indicate injection wells and the negative numbers of days indicate production wells.

The gas saturation maps $S$ (second column in Figure~\ref{fig:data}) and pressure maps $P$ (third column in Figure~\ref{fig:data}) are simulated with the state-of-the-art full-physics numerical simulator ECLIPSE \cite{eclipse}. Each simulation takes around 2 minutes on an Intel Core i7-4790 CPU.

Finally, we concatenate a sets of $\mathbf y$, $S$, and $P$ with a set of $\mathbf y'$, $S'$, and $P'$ (permuted) to become one data tuple as shown in Figure~\ref{fig:data}. The $\mathbf y$, $S$, $P$, $\mathbf y'$, $S'$, and $P'$ within each tuple are always correlated to the same permeability map $\mathbf k$ but different well location maps $\mathbf y$. 

\subsection{Training, validation, and test sets} 
For the training and the validation set, a total of 42 different permeability map realizations $\mathbf k$ (under the same Gaussian distribution) and 50 different well locations maps $\mathbf y$ were simulated in $42 \times 50 = 2100$ runs. Each simulation lasts for 1,000 days and we took 4 time snapshots of the output, generating $2100 \times 4 = 8400$ sets of $\mathbf y$, $S$, and $P$ data and $8400$ data tuples. Among the $8400$ data tuples, we eliminated $16$ tuples with extremely large pressure data, which are caused by numerical instability, and got a total of $8384$ data tuples for training and validation. The training/validation set split is 7600/784 (roughly 10:1). 

For the test set, we considered a total of 100 different permeability map realizations $\mathbf k$ (under the same Gaussian distribution as before) and 5 different well locations maps $\mathbf y$, simulated in $100 \times 5 = 500$ runs. The tuple preparation process is similar to the training and the validation set. Since our test results is compared with a conventional MC method (more details in Section 4.1), the required number of permeability map scenarios are significantly higher than in the training set. In fact, rigorous conventional MC uncertainty quantification often requires more than 5,000 of different permeability map scenarios, which is computationally intractable in the scope of this project. The drawbacks of using a small number of permeability maps in the conventional MC is further discussed in Section 6.

\section{Methods}\label{sec::methods}
Our project is built up in two stages. In the first stage, we focus on replacing conventional PDE solver with neural network. For a deterministic system, the saturation map is predicted for a given pair of permeability field and well location using U-Net. Basically, it is supervised learning in the relationship of $\mathbf k,\mathbf y$ and $S$ governed by equation~\ref{eq:3}. Satisfactory results have be shown in section 5 of previous milestone report.

The second stage is to predict the saturation map for a given well location but without any prescribed knowledge of permeability. This is a more challenging and practical problem in CCS project. The detailed geological information of the field, e.g. the permeability, is always incomplete and uncertain. This project will be useful in the evaluation of the optimal well location even in the absence of the permeability information. In a real CCS project, data of the saturation map and pressure map will be collected from test well locations (usually not the optimal sites) in a pilot test. Geological information is then interpreted from the pilot test data to build up stochastic models. Then simulation is conducted to evaluate potential optimal well controls based on stochastic models. Our project aims to interpret the permeability from pilot test saturation and pressure map in existing test wells and predict the saturation maps in other potential well locations using the same deep learning framework. In other words, the network would be able to predict the saturation map in an arbitrary well location input different from the training sets. Efficient quantification can be conducted in the evaluation of optimal well controls.

\subsection{Conventional MC Review}
We briefly review the conventional MC approach here and employ the MC solution as a baseline to compare with in the later experiment section. For each different well locations $\{\mathbf y_1,\mathbf y_2,\cdots,\mathbf y_N\}$, one draw a set of permeability sample maps $\{\mathbf k_1(\mathbf x),\mathbf k_2(\mathbf x),\cdots,\mathbf k_M(\mathbf x)\}$ following the same Gaussian distribution. For these $N\times M$ pair $(\mathbf y_i,\mathbf k_j)$, the following deterministic system is solved:
\begin{equation}\label{eq:3}
\begin{aligned}
&\phi\partial_t S+\nabla\cdot \mathbf u =0\\
&\mathbf u(\mathbf x) = -k_j(\mathbf x)\nabla \Phi(\mathbf x), \mathbf x\in\mathcal D,\\
&\nabla\cdot \mathbf u(\mathbf x) = f(\mathbf x),\mathbf x\in\mathcal D,\\
&\mathbf u(\mathbf x) \cdot \hat{\mathbf n}(\mathbf x)=0 , \mathbf x\in\partial D,\\
&f(\mathbf x) = \left\{
    \begin{aligned}
    &Q, &\mbox{if} \ \mathbf x= \mathbf y_i,\\
    &0, &\mbox{otherwise}.
    \end{aligned}
    \right.
\end{aligned}
\end{equation} 
The solution is recorded as $S_{i,j}$ for different injection days: $\{S_{i,j}^{t_1},S_{i,j}^{t_2},\cdots\}$. Notice that the saturation map under each well control is completely independent of each other. For a chosen well control $\mathbf y_i$ in a fixed injection day $t_n$, one can investigate quantities of interest from the sample set $\{S_{i,1}^{t_n},S_{i,2}^{t_n},\cdots,S_{i,M}^{t_n}\}$. For example, the mean saturation map $\bar S_i^{t_n}$ can be approximated by
\begin{equation}
\bar S_i^{t_n} \approx \frac{1}{M}\sum_{j=1}^MS_{i,j}^{t_n}.
\end{equation}
Law of large numbers guarantees the convergence of the above approximation with a rate of $O(M^{-1/2})$.

\subsection{Basic U-Net}
The first step of our project is to do supervised learning using basic U-Net~\cite{ronneberger2015u, zhu2018bayesian}. In~\cite{zhu2018bayesian}, a more appropriate network architecture, modified from the conventional U-Net~\cite{ronneberger2015u}, is applied to single-phase flow problem. With the similarity of single-phase flow and our CO$_2$ injection, the U-Net architecture from~\cite{zhu2018bayesian} is adapted here. We simulate input images of $\{(\mathbf y_i,\mathbf k_i)\}$ and the corresponding output target images $\{S_{i}\}$ by equation (3). The architecture is illustrated in the following Figure~\ref{fig:unet}:
\begin{figure}[!htp]
    \centering
    \includegraphics[scale=0.23]{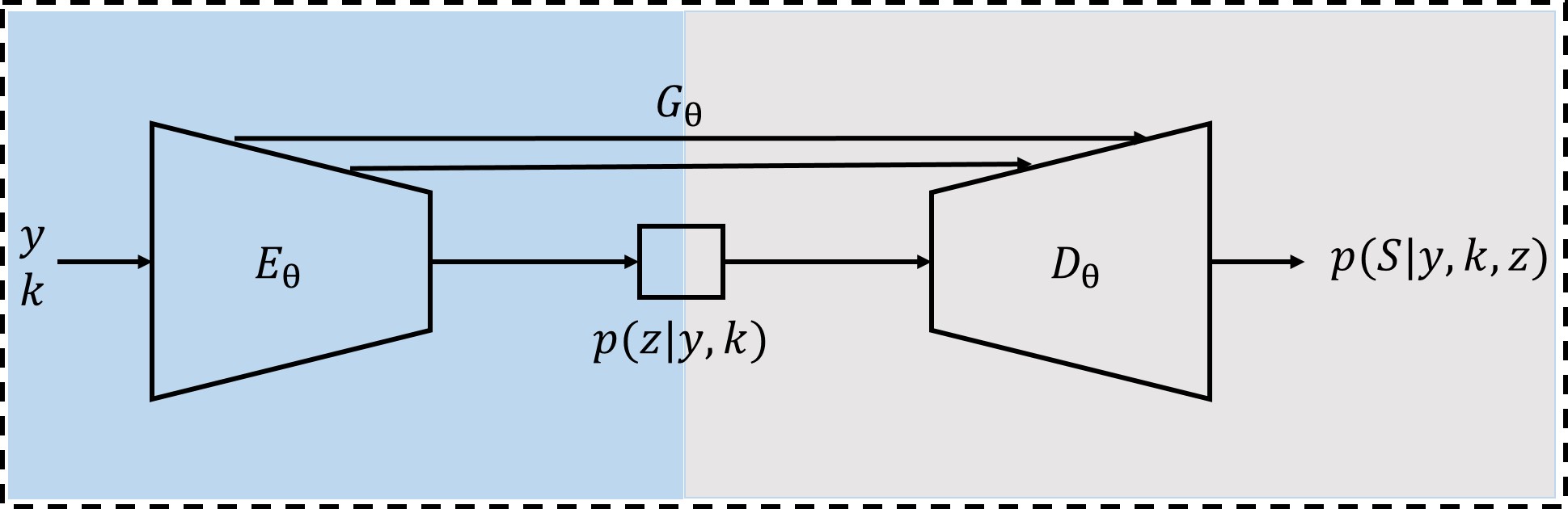}
    \caption{U-Net architecture}
    \label{fig:unet}
\end{figure}
Here $z$ is the latent variable. The following regularized MSE loss is employed as the training loss function:
\begin{equation}\label{eq:5}
    L(\mathbf f_\theta(\mathbf y,\mathbf k),S) = \frac{1}{n}\sum_{i=1}^n[\mathbf f_\theta(\mathbf y_i, \mathbf k_i)-S_i]^2+\lambda\|\theta\|_2,
\end{equation}
where $\mathbf f_\theta$ is the neural net depending on the parameter $\theta$ and $\lambda$ is the regularization strength. We refer to Fig 1 in~\cite{ronneberger2015u} for the detailed architecture of $\mathbf f_\theta$.

The performance of the prediction is quantified through the mean relative error, i.e.,
\begin{equation}\label{eq:error}
    e = \frac{1}{n}\sum_{i=1}^n\frac{\|\mathbf{f}_\theta(\mathbf y_i,\mathbf k_i)-S_i\|_{2}^{2}}{\|S_i\|_{2}^{2}},
\end{equation}
where $n$ is the total number of data points. Notice that everything is deterministic and explicitly forward (i.e., permeability information is needed as an input) at this point. The achievement of basic U-Net is to find a neural network that can replace the PDE solver for eqn~(\ref{eq:3}) efficiently. The results of Basic U-Net are presented in our milestone reports.

\subsection{Variational U-Net}
Motivated by \cite{esser2018variational}, where a conditional U-Net for shape-guided image generation is combined with a VAE for appearance conditioned on the output, we find an analogy in \cite{esser2018variational} and our problem:
\begin{equation}
\begin{aligned}
&\mbox{image } \longleftrightarrow \mbox{saturation map,}\\
&\mbox{shape } \longleftrightarrow \mbox{well location,}\\
&\mbox{appearance } \longleftrightarrow \mbox{permeability.}
\end{aligned}
\end{equation}

 Since the well location and the permeability capture all variations of saturation, the architecture proposed in \cite{esser2018variational} should be able to provide a saturation map generator controlled by the well location and permeability. We refer to Figure 2 and equation (3) in \cite{esser2018variational} for the architecture and loss function construction. However, the test results are not satisfactory. As we went through the loss function construction and the assumptions, we realized the potential cause of the problem. One key term in the ELBO of the loss function construction is the prior of appearance conditioned on shape, which captures potential interrelations between shape and appearance. However, there is no correlation between well location and the permeability field. These two variables are completely independent in our problem. Another issue is the assumption of Gaussian distribution of appearance conditioned on image and shape, which could be satisfied only if we change our previous training dataset.

\begin{figure}[!htp]
  \centering
  \includegraphics[scale = 0.23]{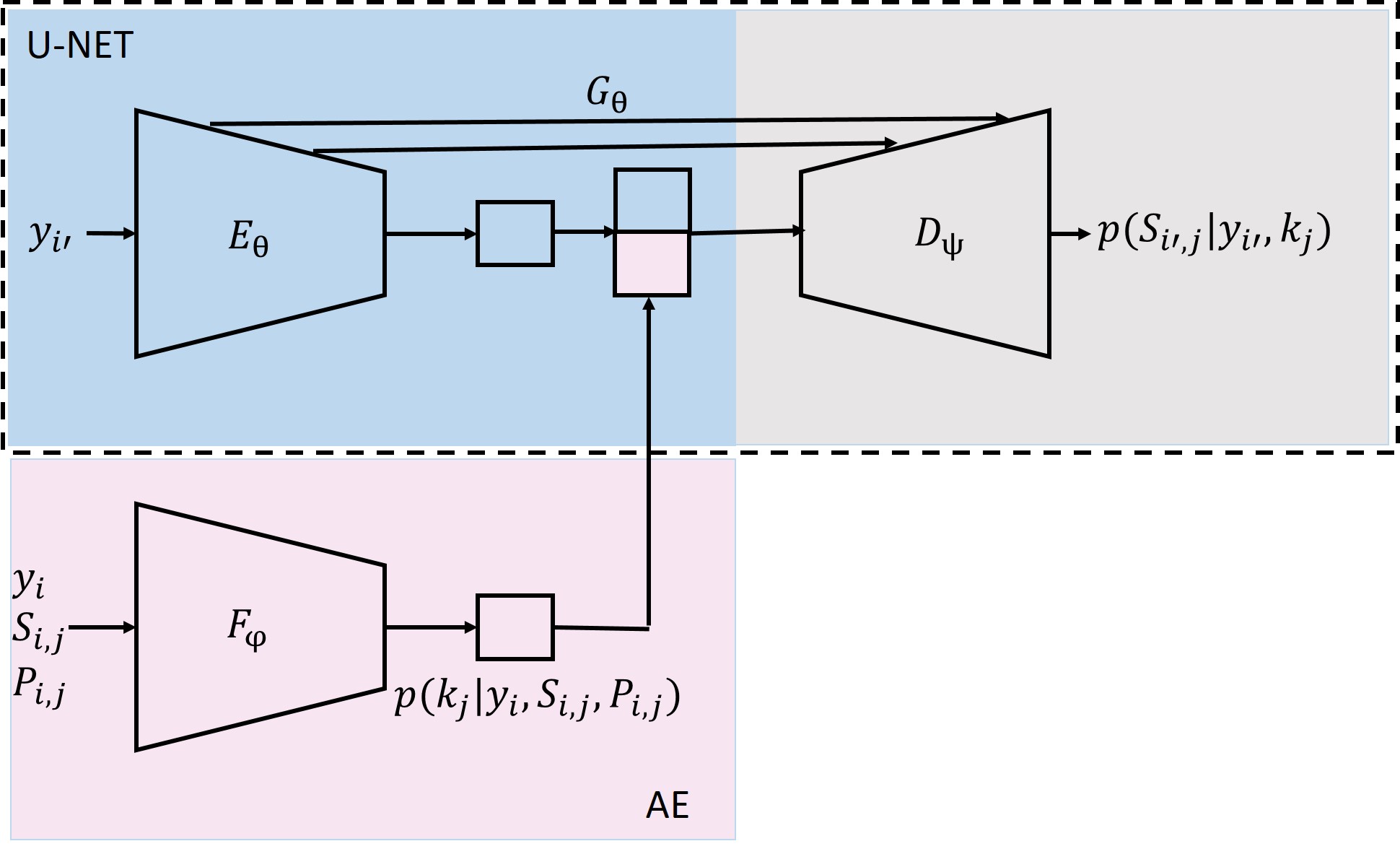}
  \caption{(Modified) variational U-Net architecture}
  \label{fig:vunet3}
\end{figure}

Figure \ref{fig:vunet3} shows our idea of the modified architecture. For a pair of input data set $(\mathbf y_i, S_{i,j},P_{i,j})$, which is saturation and pressure maps obtained from simulations in the existing old well locations $\mathbf y_i$, autoencoder will learn the latent permeability $\mathbf k_j$ conditioned on  $(\mathbf y_i, S_{i,j},P_{i,j})$. Here we involve the pressure map $P_{i,j}$, which can be conveniently obtained during the process, to better infer the underlying permeability. Then a U-Net will infer the new well location $\mathbf y_{i'}$ combined with the autoencoder so that the decoder in U-Net can generate a prediction of the new saturation map based on the new well location $\mathbf y_{i'}$ and the latent permeability $\mathbf k_j$.

The loss function is modified as below:
\begin{equation}\label{eq:4}
\begin{aligned}
    L(\mathbf f_{\theta,\phi,\psi}(\mathbf y_{i'}),S_{i',j}) 
    & = \frac{1}{n}\sum_{i=1}^n[\mathbf f_\phi(\mathbf y_i)-S_{i,j}]^2 \\
    & +\lambda_\theta\|\theta\|_2+\lambda_\phi\|\phi\|_2+\lambda_\psi\|\psi\|_2,
\end{aligned}
\end{equation}
where $\lambda_{\theta}$,$\lambda_{\phi}$ and $\lambda_{\psi}$ are the regularization strength for $\theta$, $\phi$ and $\psi$ respectively. The performance of the predict, again is evaluated with the average relative error defined in Eq.~\ref{eq:error}. Notice that there the U-Net is combined with an AE not a VAE (as in \cite{esser2018variational}), and hence the loss function is defined as MSE loss here instead of the perceptual loss in \cite{esser2018variational}. The mean relative error of saturation (and similarly pressure) is defined the same as in eqn (6). The detailed implementation of $E_{\theta}$, $D_{\psi}$ and $F_{\phi}$ can be found in Appendix.

Compared with conventional MC approach, the benefits of the proposed method in computation efficiency come from two parts. First, the AE part substitute conventional PDE solver. Notice that the permeability information is interpreted from AE and not needed as an input anymore, which solves an inverse problem of importance already. Second, new well control is directly combined with underlying permeability in the network. The prediction will be saturation under a completely new well control that has never been solved by PDE or trained in the network. This enables freedom of well control without any additional computational costs (once training is complete). In the test process, one only needs to input saturation and pressure map from existing old well location, then the trained network will infer all the uncertainties of the underlying permeability and propagate the uncertainties to the saturation map in the new selected well locations. To achieve the same goal, traditional MC has no smart way but run expensive simulations again under different well locations in each permeability sample space. Hence, V-UNET enables us to get the uncertainty information in a fraction of time.

Quantity of statistic interests can be computed in the same fashion as conventional MC approach. For example, the mean saturation map $\bar S_{i'}^{t_n}$ in the new well location $\mathbf y_{i'}$ at injection day $t_n$ can be approximated by
\begin{equation}
    \bar S_{i'}^{t_n} \approx \frac{1}{M}\sum_{j=1}^M S_{i',j}^{t_n}.
\end{equation}



\section{Experiments}

The implementation of the variational U-Net following \cite{esser2018variational} is made with modifications on the model architecture. The model is trained with 7600 data points (train set). The validation set (dev. set) contains 784 data points. We test the results for UQ task on a separate data set of 500 data points (test set).

The quantities of interest in this context are the pressure and CO$_2$ saturation fields of the reservoir. As is mentioned in Section~\ref{sec::methods}, permeability fields are only used for the full-physics numerical simulator to generate quantities of interest, and is not an input of our V-UNet model. We use mean relative error defined by Eq.~\ref{eq:error} as a metric of the performance. Note that the test error (shown in Table~\ref{tab:error} and \ref{tab:ablation}) is the average of the mean relative error for both pressure and saturation.

We first run the experiments to predict the pressure and saturation fields for a single geological realization. A variety of model architectures and hyper-parameters are tested and compared. We also conducted the ablation study on the resNet layers between conv layers and at the low-dimension (see Appendix for detail) to see their impacts on the model. Finally, we compare our V-UNet based MC model with the conventional MC (baseline model) to find the overall speedup achieved and the accuracy loss. The results and discussion for each part are shown subsequently.

\subsection*{Single realization}

We first run the experiments to predict the pressure and saturation fields for a single geological realization. A variety of model architectures and hyper-parameters are tested and compared, which are shown in Table~\ref{tab:error}. Note that the `normal' architecture in the table refers to the architecture demonstrated in Figure~\ref{fig:vunet3}. The `reversed' architecture refers to the model that $E_{\theta}$ and $F_{\phi}$ in Figure~\ref{fig:vunet3} are in reversed order (skip connections are established between $F_{\phi}$ and $D_{\psi}$). Models with varying learning rates (lr) and training epochs (ep) are also listed in the table.

We point out that one of the model settings does not give a converging loss value, so the associated test error is not shown. We can tell that the `normal' architecture with a learning rate of 2e-4 and epoch of 10 gives the lowest dev. loss, and the test error associated is 12.7\%, which is a reasonably good accuracy for surrogate models in general. We will use this model setting as the `default' setting for the following ablation study.

\begin{table}[!htp]
\centering
\begin{tabular}{ l@{\hskip .15in} l@{\hskip .05in} l@{\hskip .05in} l } 
 \hline
 Model \& Param       & Train loss & Dev. loss & Test err. \\ 
 \hline
 Normal (lr 5e-4, ep 15)    &  263.5    & 359.1    & 13.2\%         \\
 Normal (lr 1e-3, ep 15)    &  2335     & 2327     & --         \\
 Reversed (lr 5e-4, ep 15)  &  358.0    & 534.9    & 18.6\%         \\ 
 Normal (lr 2e-4, ep 15)    &  212.3    & 739.4    & 22.4\%        \\
 Normal (lr 2e-4, ep 10)    &  257.5    & \textbf{283.1}    & \textbf{12.7\%}     \\ 
 
 \hline
\end{tabular}
\caption{Model architecture and hyper-parameter tuning}
\label{tab:error}
\end{table}

\begin{figure}[!htp]
    \centering
    \includegraphics[width=0.5\textwidth]{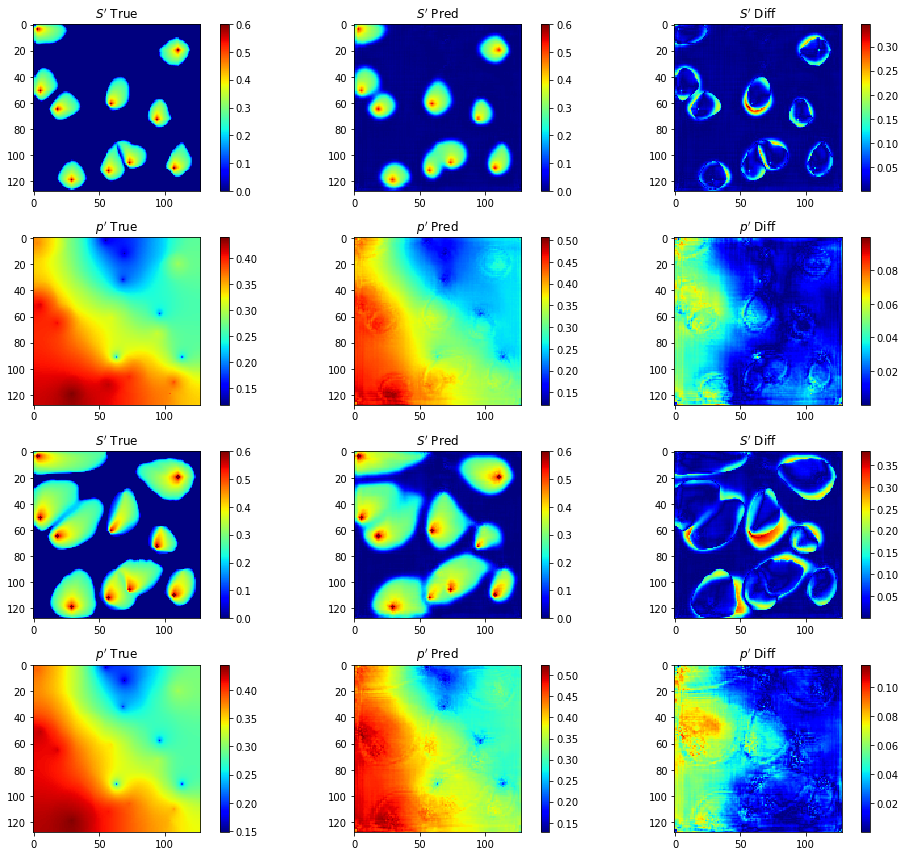}
    \caption{Test simulation results with V-UNet for prediction of single realization (Row 1 and 2: saturation and pressure at 250 days. Row 3 and 4 saturation and pressure at 750 days. Column 1: True, Column 2: predicted, Column 3: difference map of Col 1 and Col 2.)}
    \label{fig:vunet_test}
\end{figure}

A set of V-UNet results (with last set model architecture and hyper-parameters in Table~\ref{tab:error}) are shown in Figure~\ref{fig:vunet_test}. We show true and predicted saturation and pressure fields at two different time steps in the operation period. Row 1 and 2 show the saturation and pressure fields at 250 days, and row 3 and 4 show the saturation and pressure fields at 750 days. Column 1 for all rows shows the true solution from the data set, Column 2 shows the fields predicted by the V-UNet model, and Column 3 shows the difference map between Column 1 and 2.

As we can tell from the figures, the V-UNet predictions are in general in good agreement with the true fields, which clearly indicates the capability of V-UNet in capturing the saturation and pressure responses under a different set of well locations. The errors for saturation mostly comes from the edges of the CO$_2$ plumes, which could potentially be an indication of issues in predicting the underlying permeability fields or that the loss function is not sensitive to the location of the edge of the plume. The pressure field predictions are visually in close match with true solutions. However, we do see some `non-smoothness' in the predictions (more obvious in the difference maps), which indicates a violation of underlying PDEs. A physically-informed loss function (as is discussed in \cite{raissi2019physics}) is a potential remedy of the issue.

\subsection*{Ablation analysis}

We perform an ablation analysis on the default V-UNet model, which comes from the last entry of Table~\ref{tab:error}. First, the resNet layers between each of the down-sampling conv layers in the encoder and between each of the up-sampling conv layers in the decoder are removed from the default model. Next, the resNet layers in the low dimension (the last three resNet layers in the encoder, and first three resNet layers in the decoder) are removed from the default model. The results of the ablation study is shown in Table~\ref{tab:ablation}. The ablated models give a higher dev. loss and test error compared to the default model, which indicates that both components are essential for the performance of the V-UNet. 

\begin{table}[!htp]
\centering
\begin{tabular}{ l@{\hskip .25in} l@{\hskip .25in} l } 
 \hline
 Model                      & Dev. Loss & Test Error \\ 
 \hline
 Default                    & 283.1     & \textbf{12.7\%} \\ 
 No resNet between layers   & 286.6     & 12.8\% \\ 
 No resNet in low dim       & 373.4     & 13.3\% \\ 
 \hline
\end{tabular}
\caption{Ablation study}
\label{tab:ablation}
\end{table}

\subsection*{UQ task compared to conventional MC}

The V-UNet model is used as a surrogate model in the Monte-Carlo (MC) framework for the uncertainty quantification (UQ) task. The resulting V-UNet-based MC is compared to conventional MC (baseline model) to get a speedup versus accuracy loss. As a reminder, the goal of this task is that given the uncertainty information (e.g., simulation results with 100 different permeability realization) under a set of well locations, and we want to get the uncertainty for quantities of interest (e.g., pressure and saturation fields) under a different set of well locations. The conventional MC will have no other way but to re-run the simulations for new well locations, which is considerably time-consuming. However, the V-UNet will be able to infer the uncertainty results (with some accuracy loss) in a fraction of time.

The comparison of UQ results are shown in Table~\ref{tab:mc}. As is mentioned earlier, the test data set contains 500 data points, which are 100 realizations for five of the well locations each. The average time in Table~\ref{tab:mc} refers to the time required for simulation of all 100 realizations corresponding to one well locations and averaged over five well locations. The average error for pressure (P) and saturation (S) are calculated analogously.

Note that the baseline model (conventional MC) are used as reference solution, which means the errors on pressure and saturation fields are defined as 0.0\%. The errors for V-UNet MC are calculated based on the difference to those of conventional MC. The error for pressure is 4.6\% and that of the saturation is 21.84\%, both of which are considered reasonable for surrogate models. The average time required for conventional MC is 12,200 seconds, while the that for V-UNet MC is 14 seconds, which corresponds to a \textbf{871} times speedup. We want to point out that a typical MC procedure will require at least 5000 simulation runs for a setting (we run 100 each due to the limitation on project time and computation resources), and running V-UNet model currently uses only a fraction of GPU memory. The speedup will be more significant in a realistic MC setting.

\begin{table}[!htp]
\centering
\begin{tabular}{ l@{\hskip .1in} l@{\hskip .1in} l@{\hskip .1in} l } 
 \hline
 Model          & Avg. Time(s)      & Avg. Err. (p)   & Avg. Err. (S)   \\ 
 \hline
 Conven. MC     & 12,200            & 0.0\%             & 0.0\%         \\ 
 V-UNet MC      & \textbf{14}       & 4.60\%           & 21.84\%        \\ 
 \hline
\end{tabular}
\caption{Comparison to traditional MC}
\label{tab:mc}
\end{table}

\begin{figure}[!htp]
    \centering
    \includegraphics[width=0.5\textwidth]{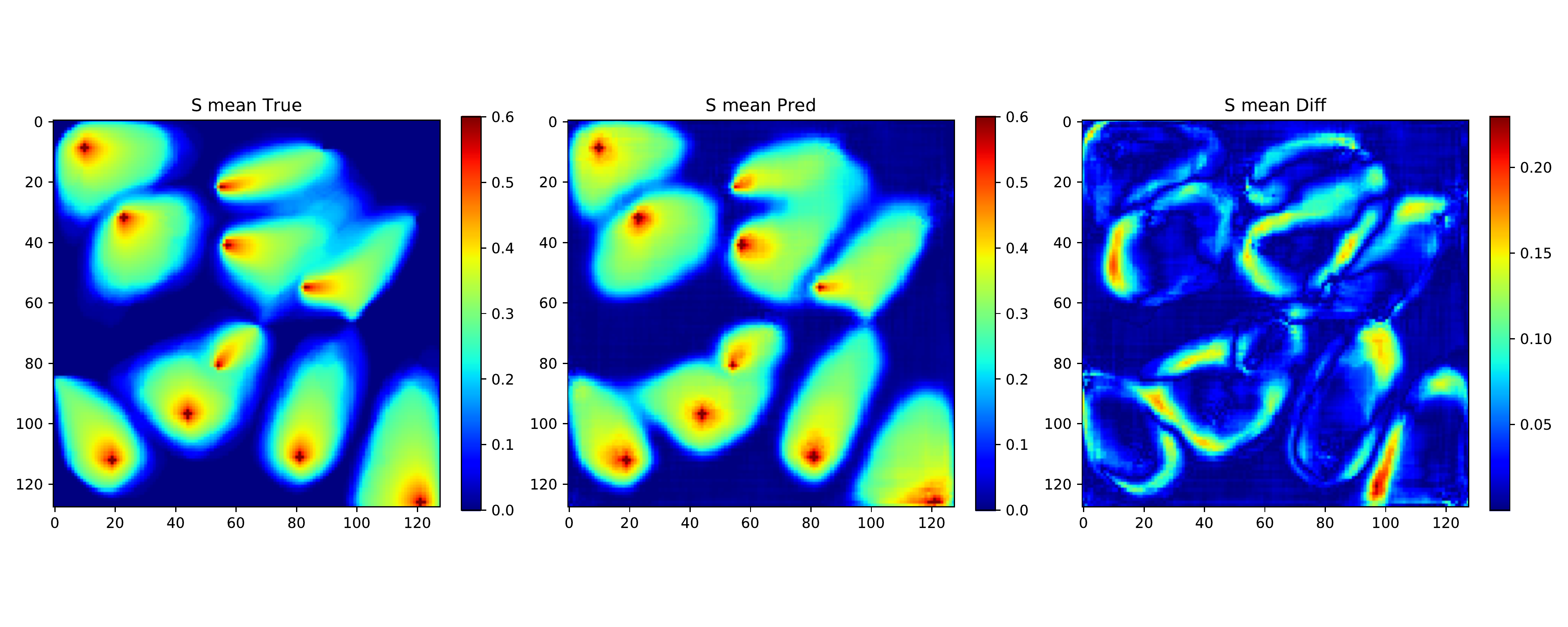}
    \caption{Mean saturation (Column 1: True, Column 2: predicted, Column 3: difference map of Col 1 and 2.)}
    \label{fig:vunet_mean_s}
\end{figure}

\begin{figure}[!htp]
    \centering
    \includegraphics[width=0.5\textwidth]{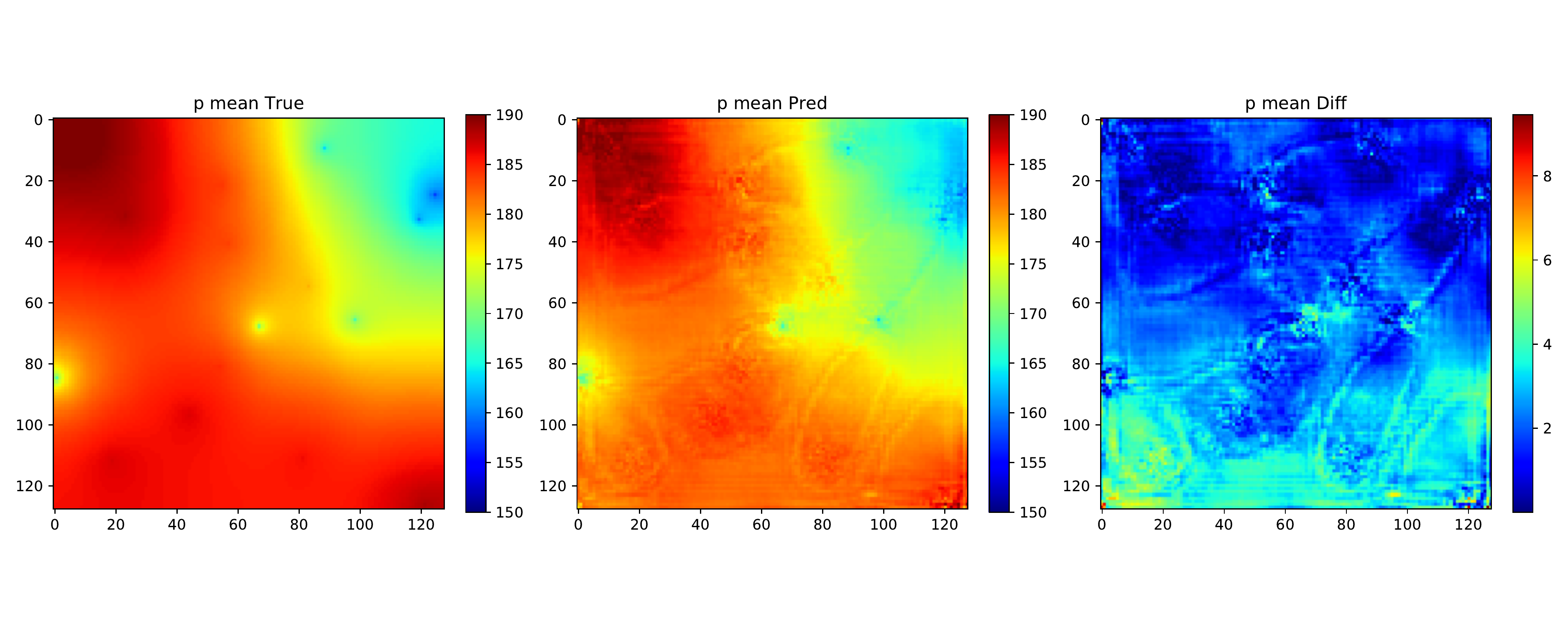}
    \caption{Mean pressure (Column 1: True, Column 2: predicted, Column 3: difference map of Col 1 and 2.)}
    \label{fig:vunet_mean_p}
\end{figure}

The mean saturation and pressure field for one of the five test well locations are shown in Figure~\ref{fig:vunet_mean_s} and \ref{fig:vunet_mean_p}. In both sets of figures, the first sub-figure shows the true solution, the second sub-figure shows the predicted solution, while the third sub-figure shows the difference map of sub-figure 1 and 2. For both sets of figures, the fields predicted by V-UNet are in close agreement with the true solution, which demonstrate the possibility of using V-UNet as a surrogate model for fast uncertainty quantification of reservoir simulation with varying well locations.

We also want to point out that the V-UNet-based MC workflow significantly under predicts the variance of both pressure and saturation. Thus the results are not shown here. The issue is still to be investigated in the future work.

\section{Conclusion and future work}
In this project, we explore the use of a variational U-NET to construct a surrogate model for a uncertain system governed by stochastic PDEs. The surrogate model is incorporated into the MC framework where quantities of statistic interest can be calculated. The image-to-image regression achieves satisfactory results in terms of predictive performance and uncertainty modeling. The main contribution of this work beyond previous studies is that we established a ``control-guided" reservoir simulation, allowing variation of well control without any additional computational costs in an equation-free framework. In a MC-based uncertainty quantification framework, the proposed method accelerates computational speed significantly (compared with conventional MC) without much sacrifice of accuracy. CO$_2$ mean saturation can be predicted in new well locations with the overall plume feature captured, providing valuable evaluation for a large-scale CCS project.

There are several directions to improve and extend in our work. First of all, the MC sampling of 100 realizations is far away from statistically enough. Limited by the computational power, 100 sample of permeability does not provide enough uncertainty information in the geological space. This could explain our observations that 1. permeability is under represented compared to well location in our surrogate model; 2. the mean saturation is well-captured but the variance and higher statistic moments are not consistent with reference solution. In the future work, we need larger dataset to do more accurate uncertainty quantification (approximately 5000 permeability sampling). Alternatively, we can design a sampling in the low dimensional space incorporated in the network to give richer statistic variations as well as more speedup. A bayesian framework with corresponding perceptual loss needs to be considered as in the original V-NET model. In the regression part, physical-informed techniques~\cite{raissi2019physics} can be added to improve prediction accuracy on pressure. Last but not least, we plan to extend our work to time-dependent problem, where Recurrent Neural Network (RNN) could be incorporated in the current framework.


\section*{Acknowledgements}

The experiments are conducted on the cluster of Stanford Center for Computational Earth \& Environmental Sciences (CEES). The provided computational resources are greatly appreciated.

\section*{Honor coder clarification}
The Github source code corresponding to \cite{esser2018variational} (\textit{https://github.com/CompVis/vunet}) is referred to during our implementation. We made significant modification with respect to the original code and wrote our version of V-UNet code from scratch (see source code in supplement material for detail).

\section*{Appendix}

Here we show the detailed architecture of encoder and decoder models used in the V-UNet framework. The architecture of encoder (same layout but different model parameters for both $E_{\theta}$ and $F_{\phi}$) is demonstrated in Table~\ref{tab::enc}. Note that $N_x = N_y = 128$ denotes the number of $x-y$ pixels in the input image. The variable $n_v = 1$ for $E_{\theta}$ since well location is the only type of input, while $n_v = 3$ for $F_{\phi}$ since we have saturation pressure and well location input at the same time.

The network-in-network layer (nin lyr) is simply a 1x1 convolutional layer with 8 filters. Both the `In res' and the `Low res' blocks denote the residual block with a structure of Conv2D-batchNorm-ReLU-Conv2D-batchNorm, where Conv2D stands for two-dimensional convolutional layer, batchNorm denotes batch normalization layer, and ReLU designates rectified linear unit. `In res' means the residual block in between each convolutional blocks, which is distinguished from `Low res' (i.e., the residual block with feature map in low dimensional space). The `Conv blk' refers to the convolutional block with a structure of conv2D-bathcNorm-ReLU.

Each of the outputs of the `In res' layer for $E_{\theta}$ will go through the skip connection of the U-Net and serve as the additional input (i.e., `skpInput') of the decoder model. The final outputs of both the encoder (with shape of $8\times8\times128$ each), will be concatenated (shape of $8\times8\times256$) and serve as the input of the decoder.

\begin{table}[!htb]
\centering
\begin{tabular}{l@{\hskip .05in} l@{\hskip .05in} l}
    \hline
    \textbf{Layer}  & \textbf{Filters}              &\textbf{Output size} \\
    \hline
    Input           &                               & ($N_{x}$, $N_{y}$, $n_v$)  \\
        \\[-1em]
    Nin lyr  & 8 of (1,1,$n_v$), stride 1           & ($N_{x}$, $N_{y}$, $8$)  \\
        \\[-1em]
    In res  & 8 of (3,3,8), stride 1               & ($N_{x}$, $N_{y}$, $8$)  \\
        \\[-1em]
    Conv blk & 16 of (3,3,8), stride 2              & ($N_{x}/2$, $N_{y}/2$, $16$)  \\
        \\[-1em]
    In res  & 16 of (3,3,16), stride 1             & ($N_{x}/2$, $N_{y}/2$, $16$)  \\
        \\[-1em]
    Conv blk & 32 of (3,3,16), stride 2             & ($N_{x}/4$, $N_{y}/4$, $32$)  \\
        \\[-1em]
    In res  & 32 of (3,3,32), stride 1             & ($N_{x}/4$, $N_{y}/4$, $32$)  \\
        \\[-1em]
    Conv blk & 64 of (3,3,32), stride 2             & ($N_{x}/8$, $N_{y}/8$, $64$)  \\
        \\[-1em]
    In res  & 64 of (3,3,64), stride 1             & ($N_{x}/8$, $N_{y}/8$, $64$)  \\
        \\[-1em]
    Conv blk & 128 of (3,3,64), stride 2            & ($N_{x}/16$, $N_{y}/16$, $128$)  \\
        \\[-1em]
    Low res  & 128 of (3,3,128), stride 1           & ($N_{x}/16$, $N_{y}/16$, $128$)  \\
        \\[-1em]
    Low res  & 128 of (3,3,128), stride 1           & ($N_{x}/16$, $N_{y}/16$, $128$)  \\
        \\[-1em]
    Low res  & 128 of (3,3,128), stride 1           & ($N_{x}/16$, $N_{y}/16$, $128$)  \\
    \hline  
\end{tabular}
\caption{Network architecture for encoder ($E_{\theta}$, $F_{\phi}$)}
\label{tab::enc}
\end{table}

The detailed structure of the decoder (i.e, $D_{\psi}$) is shown in Table~\ref{tab::dec}. The `low res' and `in res' are defined similarly as those in the encoder. The `upConv' stands for the upconvolutional block with a sequential stack of a two-dimensional unpooling layer, a reflection padding layer, a two-dimensional convolutional layer, a batch normalization layer, and a ReLU, which in the end expand the spatial dimension by two. The output of the upConv block is concatenated with the input from skip connection (`skpInput') with the exact same dimension before input to the next residual block. We have a conv2D layer in the end to convert the feature map to saturation and pressure fields.

\begin{table}[!htb]
\centering
\begin{tabular}{l@{\hskip .05in} l@{\hskip .05in} l}
    \hline
    \textbf{Layer}  & \textbf{Filters}                       &\textbf{Output size} \\
    \hline
    Input           &                                   & ($N_{x}/16$, $N_{y}/16$, $256$)  \\
        \\[-1em]
    Low res & 256 of (3,3,256), stride 1                & ($N_{x}/16$, $N_{y}/16$, $256$)  \\
        \\[-1em]
    Low res & 256 of (3,3,256), stride 1                & ($N_{x}/16$, $N_{y}/16$, $256$)  \\
        \\[-1em]
    Low res & 256 of (3,3,256), stride 1                & ($N_{x}/16$, $N_{y}/16$, $256$)  \\
        \\[-1em]
    upConv  & 64 of (3,3,256), stride 2                 & ($N_{x}/8$, $N_{y}/8$, $64$)  \\
        \\[-1em]
    skpInput& concat. upConv output                     & ($N_{x}/8$, $N_{y}/8$, $64$)  \\
        \\[-1em]
    In res  & 128 of (3,3,128), stride 1                & ($N_{x}/8$, $N_{y}/8$, $128$)  \\
        \\[-1em]
    upConv  & 32 of (3,3,128), stride 2                 & ($N_{x}/4$, $N_{y}/4$, $32$)  \\
        \\[-1em]
    skpInput& concat. upConv output                     & ($N_{x}/4$, $N_{y}/4$, $32$)  \\
        \\[-1em]
    In res  & 64 of (3,3,64), stride 1                  & ($N_{x}/4$, $N_{y}/4$, $64$)  \\
        \\[-1em]
    upConv  & 16 of (3,3,64), stride 2                  & ($N_{x}/2$, $N_{y}/2$, $16$)  \\
        \\[-1em]
    skpInput& concat. upConv output                     & ($N_{x}/2$, $N_{y}/2$, $16$)  \\
        \\[-1em]
    In res  & 32 of (3,3,32), stride 1                  & ($N_{x}/2$, $N_{y}/2$, $32$)  \\
        \\[-1em]
    upConv  & 8 of (3,3,32), stride 2                   & ($N_{x}$, $N_{y}$, $8$)  \\
        \\[-1em]
    skpInput& concat. upConv output                     & ($N_{x}$, $N_{y}$, $8$)  \\
        \\[-1em]
    In res  & 16 of (3,3,16), stride 1                  & ($N_{x}$, $N_{y}$, $16$)  \\
        \\[-1em]
    Conv lyr  & 2 of (3,3,16), stride 1                 & ($N_{x}$, $N_{y}$, $2$)  \\

    \hline  
\end{tabular}
\caption{Network architecture for decoder ($D_{\psi}$)}
\label{tab::dec}
\end{table}

{\small
\bibliographystyle{ieee}
\bibliography{main}

\begin{thebibliography}{10}\itemsep=-1pt

\bibitem{babuvska2007stochastic}
I.~Babu{\v{s}}ka, F.~Nobile, and R.~Tempone.
\newblock A stochastic collocation method for elliptic partial differential
  equations with random input data.
\newblock {\em SIAM Journal on Numerical Analysis}, 45(3):1005--1034, 2007.

\bibitem{baldi2014searching}
P.~Baldi, P.~Sadowski, and D.~Whiteson.
\newblock Searching for exotic particles in high-energy physics with deep
  learning.
\newblock {\em Nature communications}, 5:4308, 2014.

\bibitem{sgems}
A.~Boucher, J.~Wu, and N.~Remy.
\newblock Applied geostatistics with sgems.
\newblock 2009.

\bibitem{burhenne2011sampling}
S.~Burhenne, D.~Jacob, and G.~P. Henze.
\newblock Sampling based on sobol?sequences for monte carlo techniques applied
  to building simulations.
\newblock In {\em Building Simulation 2011: 12th Conference of International
  Building Performance Simulation Association, Sydney, Australia, Nov}, pages
  14--16, 2011.

\bibitem{chan2018machine}
S.~Chan and A.~H. Elsheikh.
\newblock A machine learning approach for efficient uncertainty quantification
  using multiscale methods.
\newblock {\em Journal of Computational Physics}, 354:493--511, 2018.

\bibitem{esser2018variational}
P.~Esser, E.~Sutter, and B.~Ommer.
\newblock A variational u-net for conditional appearance and shape generation.
\newblock In {\em Proceedings of the IEEE Conference on Computer Vision and
  Pattern Recognition}, pages 8857--8866, 2018.

\bibitem{giles2015multilevel}
M.~B. Giles.
\newblock Multilevel monte carlo methods.
\newblock {\em Acta Numerica}, 24:259--328, 2015.

\bibitem{goodfellow2014generative}
I.~Goodfellow, J.~Pouget-Abadie, M.~Mirza, B.~Xu, D.~Warde-Farley, S.~Ozair,
  A.~Courville, and Y.~Bengio.
\newblock Generative adversarial nets.
\newblock In {\em Advances in neural information processing systems}, pages
  2672--2680, 2014.

\bibitem{IPCC}
IPCC.
\newblock Climate change 2014: Synthesis report - summary chapter for
  policymakers.
\newblock page~31, 2014.

\bibitem{kalla2005use}
S.~Kalla.
\newblock Use of orthogonal arrays, quasi-monte carlo sampling and kriging
  response models for reservoir simulation with many varying factors.
\newblock 2005.

\bibitem{kingma2013auto}
D.~P. Kingma and M.~Welling.
\newblock Auto-encoding variational bayes.
\newblock {\em arXiv preprint arXiv:1312.6114}, 2013.

\bibitem{kutz2017deep}
J.~N. Kutz.
\newblock Deep learning in fluid dynamics.
\newblock {\em Journal of Fluid Mechanics}, 814:1--4, 2017.

\bibitem{lu2016improved}
D.~Lu, G.~Zhang, C.~Webster, and C.~Barbier.
\newblock An improved multilevel monte carlo method for estimating probability
  distribution functions in stochastic oil reservoir simulations.
\newblock {\em Water resources research}, 52(12):9642--9660, 2016.

\bibitem{marccais2017prospective}
J.~Mar{\c{c}}ais and J.-R. De~Dreuzy.
\newblock Prospective interest of deep learning for hydrological inference.
\newblock {\em Groundwater}, 55(5):688--692, 2017.

\bibitem{min2017deep}
S.~Min, B.~Lee, and S.~Yoon.
\newblock Deep learning in bioinformatics.
\newblock {\em Briefings in bioinformatics}, 18(5):851--869, 2017.

\bibitem{mo2018deep}
S.~Mo, Y.~Zhu, N.~Zabaras, X.~Shi, and J.~Wu.
\newblock Deep convolutional encoder-decoder networks for uncertainty
  quantification of dynamic multiphase flow in heterogeneous media.
\newblock {\em arXiv preprint arXiv:1807.00882}, 2018.

\bibitem{raissi2019physics}
M.~Raissi, P.~Perdikaris, and G.~E. Karniadakis.
\newblock Physics-informed neural networks: A deep learning framework for
  solving forward and inverse problems involving nonlinear partial differential
  equations.
\newblock {\em Journal of Computational Physics}, 378:686--707, 2019.

\bibitem{robert2013monte}
C.~Robert and G.~Casella.
\newblock {\em Monte Carlo statistical methods}.
\newblock Springer Science \& Business Media, 2013.

\bibitem{ronneberger2015u}
O.~Ronneberger, P.~Fischer, and T.~Brox.
\newblock U-net: Convolutional networks for biomedical image segmentation.
\newblock In {\em International Conference on Medical image computing and
  computer-assisted intervention}, pages 234--241. Springer, 2015.

\bibitem{eclipse}
Schlumberger.
\newblock Eclipse reference manual.
\newblock 2014.

\bibitem{tartakovsky2016method}
D.~M. Tartakovsky and P.~A. Gremaud.
\newblock Method of distributions for uncertainty quantification.
\newblock {\em Handbook of Uncertainty Quantification}, pages 1--22, 2016.

\bibitem{winter2003moment}
C.~L. Winter, D.~Tartakovsky, and A.~Guadagnini.
\newblock Moment differential equations for flow in highly heterogeneous porous
  media.
\newblock {\em Surveys in Geophysics}, 24(1):81--106, 2003.

\bibitem{xiu2010numerical}
D.~Xiu.
\newblock {\em Numerical methods for stochastic computations: a spectral method
  approach}.
\newblock Princeton university press, 2010.

\bibitem{xiu2003modeling}
D.~Xiu and G.~E. Karniadakis.
\newblock Modeling uncertainty in flow simulations via generalized polynomial
  chaos.
\newblock {\em Journal of computational physics}, 187(1):137--167, 2003.

\bibitem{zhu2018bayesian}
Y.~Zhu and N.~Zabaras.
\newblock Bayesian deep convolutional encoder--decoder networks for surrogate
  modeling and uncertainty quantification.
\newblock {\em Journal of Computational Physics}, 366:415--447, 2018.

\end{thebibliography}
}

\end{document}